\documentclass[a4paper,11pt]{amsart}
\usepackage{amssymb}

\numberwithin{equation}{section}
\newtheorem{theorem}{Theorem}[section]

\newtheorem{lem}[theorem]{Lemma}

\theoremstyle{remark}

\renewcommand{\hat}{\widehat}

\author[J.~Benameur]{Jamel Benameur}
\address{Department of Mathematics, College of Science, King Saud University\\
Riyadh 11451, Kingdom of Saudi Arabia}
\email{\sl jbenameur@ksu.edu.sa}
\author[L.~Jlali]{Lotfi Jlali}
\address{Department of Mathematics, College of Science, King Saud University\\
Riyadh 11451, Kingdom of Saudi Arabia}
\email{\sl ljlali@ksu.edu.sa}

\title[Long time decay of 3D-NSE in Lei-Lin-Gevrey spaces]
{Long time decay of 3D-NSE in Lei-Lin-Gevrey spaces}
\date{\today}

\begin{document}
\begin{abstract}
In this paper, we prove that there exists a unique global solution of $3D$ Navier-Stokes equation if $\exp(a|D|^{1/\sigma})u^0\in{\mathcal{X}}^{-1}(\mathbb R^3)$ and $\|u^0\|_{{\mathcal{X}}^{-1}}<\nu$. Moreover, we will show that $\|\exp(a|D|^{1/\sigma}) u(t)\|_{{\mathcal{X}}^{-1}}$ goes to zero if the time $t$ goes to infinity.
\end{abstract}


\subjclass[2000]{35-xx, 35Bxx, 35Lxx}
\keywords{Navier-Stokes Equations; Critical spaces; Long time decay}

\maketitle
\tableofcontents


\section{Introduction}

The $3D$ incompressible Navier-Stokes equations are given by:
$$
\left\{
  \begin{array}{lll}
     \partial_t u
 -\nu\Delta u+ u.\nabla u  & =&\;\;-\nabla p\hbox{ in } \mathbb R^+\times \mathbb R^3\\
     {\rm div}\, u &=& 0 \hbox{ in } \mathbb R^+\times \mathbb R^3\\
    u(0,x) &=&u^0(x) \;\;\hbox{ in }\mathbb R^3,
  \end{array}
\right.
\leqno{(NSE)}
$$
where $\nu>0$ is the viscosity of fluid, $u=u(t,x)=(u_1,u_2,u_3)$ and $p=p(t,x)$ denote respectively the unknown velocity and the unknown pressure of the fluid at the point $(t,x)\in \mathbb R^+\times \mathbb R^3$, and $(u.\nabla u):=u_1\partial_1 u+u_2\partial_2 u+u_3\partial_3u$, while $u^0=(u_1^o(x),u_2^o(x),u_3^o(x))$ is an initial given velocity. If $u^0$ is quite regular, the divergence free condition determines the pressure $p$.\\

The study of local existence is studied by serval researchers, Leray \cite{JL01,JL02}, Kato \cite{TK1},etc.\\
The global existence of weak solutions goes back to Leray \cite{JL02} and Hopf \cite{EH}. The global well-posedness of strong solutions for small initial data in the critical Sobolev space $\dot{H}^{\frac{1}{2}}$ is due to Fujita and Kato \cite{TK3}, also in \cite{JYC1}, Chemin has proved the case of $\dot{H}^s$, $s>\frac{1}{2}$. In \cite{TK2}, Kato has proved the case of Lebesgue space $L^3$. In \cite{HKD}, Koch and Tataru have proved the case of the space $\mathbf{B\mathbf{M\mathbf{O}^{-1}}}$ (see, also \cite{MC2,JYC2,FP}). It should be noted, in all these works, that the norms in corresponding spaces of the initial data are assumed to be very small, smaller than the viscosity $\nu$ multiplied by tiny positive constant $c$. For further results and details the reader can consult the book by Cannone \cite{MC1}. In \cite{ZL}, the authors consider a new critical space that is contained in $\mathbf{B\mathbf{M\mathbf{O}^{-1}}}$, where they show it is sufficient to assumed the norms of initial data are less than exactly the viscosity coefficient $\nu$. Then, the used space  in \cite{ZL} is the following
$$ \mathcal{X}^{-1}(\mathbb R^3) =\{f\in{\mathcal D}'(\mathbb R^3);\;\int_{\mathbb R^3} \frac {|\hat{u}(\xi)|}{|\xi|}d\xi<\infty\}$$
which is equipped with the norm $$\|f\|_{ \mathcal{X}^{-1}(\mathbb R^3)}=\int_{\mathbb R^3} \frac {|\hat{u}(\xi)|}{|\xi|} d\xi.$$
We will also use the notation, for $i=0,1$,
$$ \mathcal{X}^{i}(\mathbb R^3) =\{f\in{\mathcal D}'(\mathbb R^3);\;\int_{\mathbb R^3}|\xi|^i|\hat{u}(\xi)|d\xi<\infty\}.$$
For the small initial data, the global existence is proved in \cite{ZL}:
\begin{theorem}
(See \cite{ZL}). Let $u^0\in  \mathcal{X}^{-1}(\mathbb R^3)$, such that $\|u^0\|_{\mathcal{X}^{-1}(\mathbb R^3)}<\nu$.  Then, there is a unique  $u\in{\mathcal C}(\mathbb R^+, \mathcal{X}^{-1}(\mathbb R^3))$ such that $\Delta u\in{L^1}(\mathbb R^+, \mathcal{X}^{-1}(\mathbb R^3))$. Moreover, $\forall t\geq0$
$$
\sup_{0\leq t<\infty}\left( \|u(t)\|_{\mathcal{X}^{-1}}+(\nu-\|u^0\|_{\mathcal{X}^{-1}})\int_0^t\|\nabla u\|_{L^{\infty}}d\tau\right)\leq\|u^0\|_{\mathcal{X}^{-1}}.
$$
\end{theorem}
Moreover, in \cite{ZZ} the authors proved the local existence for the initial data and blow-up criteria if the maximal time is finite, precisely:
\begin{theorem} (See \cite{ZZ}).
Let $u^0\in  \mathcal{X}^{-1}(\mathbb R^3)$. There exists time $T$ such that the system $(NSE)$ has unique solution $ u\in{L^2}([0,T], \mathcal{X}^{0}(\mathbb R^3))$ wish also belong to
$$
{\mathcal C}([0,T], \mathcal{X}^{-1}(\mathbb R^3))\cap{L^1}([0,T], \mathcal{X}^{1}(\mathbb R^3))\cap{L^{\infty}}([0,T], \mathcal{X}^{-1}(\mathbb R^3))
$$
Let $T^*$ denote the maximal time of existence of such solution. Hence\\
If $\|u\|_{\mathcal{X}^{-1}}<\nu$, then
$$
T^*=\infty
$$
If $T^*$ is finite, then
$$
\int_0^{T^*}\|u(t)\|_{\mathcal{X}^{0}}^2=\infty.
$$
\end{theorem}
Also, the long time decay for the global solution was studied in \cite{JB}, precisely:
\begin{theorem}\label{theo1}
(See \cite{JB}) Let $u\in{\mathcal C}(\mathbb R^+, \mathcal{X}^{-1}(\mathbb R^3))$ be a global solution of $(NSE)$, then
$$
\lim\sup_{t\rightarrow\infty}\|u(t)\|_{\mathcal{X}^{-1}}=0.
$$
\end{theorem}

To prepare for announce our main results, we  need to introduce the Lei-Lin-Gevrey spaces: For $a>0$, $\sigma>1$ and $\rho\in\mathbb R $, the following spaces are defined
$$ Z^{\rho}_{a,\sigma}(\mathbb R^3) =\{f\in{\mathcal S}'(\mathbb R^3);\;\int_{\mathbb R^3}|\xi|^{\rho} e^{a|\xi|^{1/\sigma}}|\hat{f}(\xi)|d\xi<\infty\}$$
which is equipped with the norm $$\|f\|_{Z^{\rho}_{a,\sigma}(\mathbb R^3)}=\int_{\mathbb R^3}|\xi|^{\rho} e^{a|\xi|^{1/\sigma}}|\hat{f}(\xi)|d\xi.$$

Our first result is the following:

\begin{theorem}\label{theo3}
Let $u^0\in Z^{-1}_{a,\sigma}(\mathbb R^3)$, such that $\|u\|_{\mathcal{X}^{-1}(\mathbb R^3)}<\nu$.  Then, there exists a unique global solution $u\in{\mathcal C}(\mathbb R^+,Z^{-1}_{a,\sigma}(\mathbb R^3))\cap {L^1}(\mathbb R^+,Z^1_{a,\sigma}(\mathbb R^3))$ of $(NSE)$.
\end{theorem}

Our second result is as follows:
\begin{theorem}\label{theo4}
Let $u\in{\mathcal C}(\mathbb R^+,Z^{-1}_{a,\sigma}(\mathbb R^3))$ be the global solution of $(NSE)$. Then
\begin{eqnarray*}
\lim\sup_{t\rightarrow\infty}\|u(t)\|_{Z^{-1}_{a,\sigma}}=0.
\end{eqnarray*}
\end{theorem}

The paper is organized in the following way: In section $2$, we give some notations and important preliminary results. Section $3$ is devoted to prove that (NSE) is well posed in $Z^{-1}_{a,\sigma}(\mathbb R^3)$. In section $4$, we prove the existence under the condition $\|u\|_{\mathcal{X}^{-1}(\mathbb R^3)}<\nu$. Finally, in the section $5$, we state  that the norm of global solution in $Z^{-1}_{a,\sigma}(\mathbb R^3)$ goes to zero at infinity.

\section{Notations and preliminary results}
\subsection{Notations}
In this section, we collect some notations and definitions that will be used later.\\
$\bullet$ The Fourier transformation is normalized as
$$
\mathcal{F}(f)(\xi)=\widehat{f}(\xi)=\int_{\mathbb R^3}\exp(-ix.\xi)f(x)dx,\,\,\,\xi=(\xi_1,\xi_2,\xi_3)\in\mathbb R^3.
$$
$\bullet$ The inverse Fourier formula is
$$
\mathcal{F}^{-1}(g)(x)=(2\pi)^{-3}\int_{\mathbb R^3}\exp(i\xi.x)g(\xi)d\xi,\,\,\,x=(x_1,x_2,x_3)\in\mathbb R^3.
$$
$\bullet$ The convolution product of a suitable pair of function $f$ and $g$ on $\mathbb R^3$ is given by
$$
(f\ast g)(x):=\int_{\mathbb R^3}f(y)g(x-y)dy.
$$
$\bullet$ If $f=(f_1,f_2,f_3)$ and $g=(g_1,g_2,g_3)$ are two vector fields, we set
$$
f\otimes g:=(g_1f,g_2f,g_3f),
$$
and
$$
{\rm div}\,(f\otimes g):=({\rm div}\,(g_1f),{\rm div}\,(g_2f),{\rm div}\,(g_3f)).
$$
$\bullet$ Let $(B,||.||)$, be a Banach space, $1\leq p \leq\infty$ and  $T>0$. We define $L^p_T(B)$ the space of all
measurable functions $[0,t]\ni t\mapsto f(t) \in B$ such that $t\mapsto||f(t)||\in L^p([0,T])$.

\subsection{Preliminary results}
In this section, we recall some classical results and we give new technical lemmas.\\
\begin{lem}\label{lem1}
Let $f,g\in Z^{-1}_{a,\sigma}(\mathbb R^3)\cap Z^{1}_{a,\sigma}(\mathbb R^3)$. Then
$$
\|fg\|_{Z^{0}_{a,\sigma}}\leq\|f\|_{Z^{-1}_{a,\sigma}}\|g\|_{Z^{1}_{a,\sigma}}+\|f\|_{Z^{1}_{a,\sigma}}\|g\|_{Z^{-1}_{a,\sigma}}.
$$
\end{lem}
\noindent{\it Proof lemma \ref{lem1}.} We have
\begin{eqnarray*}
\|fg\|_{Z^{0}_{a,\sigma}}&=&\int_{\mathbb R^3}e^{a|\xi|^{1/\sigma}}|\widehat{fg}(\xi)|d\xi\\&\leq&\int_{\xi}e^{a|\xi|^{1/\sigma}}
\left(\int_{\eta}|\hat{f}(\xi-\eta)||\hat{g}(\eta)|d\eta\right)d\xi.
\end{eqnarray*}
Using the inequality
$e^{a|\xi|^{1/\sigma}}\leq e^{a|\xi-\eta|^{1/\sigma}}e^{a|\eta|^{1/\sigma}}$ and $1\leq \frac{|\xi-\eta|}{|\eta|}+\frac{|\eta|}{|\xi-\eta|}$.\\
We obtain
\begin{eqnarray*}
\|fg\|_{Z^{0}_{a,\sigma}}&=&\int_{\xi}(\int_{\eta}|\xi-\eta|e^{a|\xi-\eta|^{1/\sigma}}|\hat{f}(\xi-\eta)|\frac{|e^{a|\eta|^{1/\sigma}}}
{|\eta|}|\hat{g}(\eta)|d\eta\\&+&\int_{\eta}\frac{e^{a|\xi-\eta|^{1/\sigma}}}{|\xi-\eta|}|\hat{f}(\xi-\eta)|\eta||e^{a|\eta|^{1/\sigma}}|\hat{g}(\eta)|d\eta
)d\xi.
\end{eqnarray*}
Put $$F_1(\xi)=|\xi|e^{a|\xi|^{1/\sigma}}|\hat{f}(\xi)|,\,\,F_2(\xi)= \frac{e^{a|\xi|^{1/\sigma}}}{|\xi|}|\hat{f}(\xi)|,\,\,
G_1(\xi)=|\xi|e^{a|\xi|^{1/\sigma}}|\hat{g}(\xi)|\,\,\,and \,\,\, G_2(\xi)= \frac{e^{a|\xi|^{1/\sigma}}}{|\xi|}|\hat{g}(\xi)|.$$
Then
\begin{eqnarray*}
\|fg\|_{Z^{0}_{a,\sigma}}&\leq&\|F_1\ast G_2\|_{L^1}+\|F_2\ast G_1\|_{L^1}\\&\leq&\|F_1\|_{L^1}\| G_2\|_{L^1}+\|F_2\|_{L^1}\|G_1\|_{L^1}
\\&\leq&\|f\|_{Z^{-1}_{a,\sigma}}\|g\|_{Z^{1}_{a,\sigma}}+|f\|_{Z^{1}_{a,\sigma}}\|g\|_{Z^{-1}_{a,\sigma}}.
\end{eqnarray*}
\hfill $\square$
\begin{lem}\label{lem2}
Let $u\in L^{\infty}_{T}(Z^{-1}_{a,\sigma}(\mathbb R^3))\cap L^{1}_{T}(Z^{1}_{a,\sigma}(\mathbb R^3))$. Then
$$
\|\int_0^t e^{\nu(t-\tau)\Delta}div(u\otimes u)d\tau\|_{Z^{-1}_{a,\sigma}}\leq2\|u\|_{L^{\infty}_{T}(Z^{-1}_{a,\sigma})}
\|u\|_{ L^{1}_{T}(Z^{1}_{a,\sigma})}.
$$
\end{lem}
\noindent{\it Proof lemma \ref{lem2}.}
\begin{eqnarray*}
\|\int_0^t e^{\nu(t-\tau)\Delta}div(u\otimes u)d\tau\|_{Z^{-1}_{a,\sigma}}&\leq&\int_0^t\|e^{\nu(t-\tau)\Delta}div(u\otimes u)\|_{Z^{-1}_{a,\sigma}}d\tau\\&\leq&\int_0^t\|e^{\nu(t-\tau)\Delta}(u\otimes u)\|_{Z^{0}_{a,\sigma}}d\tau\\&\leq&
\int_0^t\|(u\otimes u)\|_{Z^{0}_{a,\sigma}}d\tau.
\end{eqnarray*}
Using the lemma \ref{lem1}, we obtain
\begin{eqnarray*}
\|\int_0^t e^{\nu(t-\tau)\Delta}div(u\otimes u)d\tau\|_{Z^{-1}_{a,\sigma}}&\leq&2\int_0^t\|u\|_{Z^{-1}_{a,\sigma}}
\|u\|_{Z^{1}_{a,\sigma}}\\&\leq&2\|u\|_{L^{\infty}_{T}(Z^{-1}_{a,\sigma})}\|u\|_{ L^{1}_{T}(Z^{1}_{a,\sigma})}.
\end{eqnarray*}
\hfill $\square$
\begin{lem}\label{lem3}
Let $u\in L^{\infty}_{T}(Z^{-1}_{a,\sigma}(\mathbb R^3))\cap L^{1}_{T}(Z^{1}_{a,\sigma}(\mathbb R^3))$. Then
$$
\int_0^T\|\int_0^t e^{\nu(t-\tau)\Delta}div(u\otimes u)d\tau\|_{Z^{1}_{a,\sigma}}dt\leq2\|u\|_{L^{\infty}_{T}(Z^{-1}_{a,\sigma})}\|u\|_{ L^{1}_{T}(Z^{1}_{a,\sigma})}.
$$
\end{lem}
\noindent{\it Proof lemma \ref{lem3}.}
\begin{eqnarray*}
\int_0^T\|\int_0^t e^{\nu(t-\tau)\Delta}div(u\otimes u)d\tau\|_{Z^{1}_{a,\sigma}}dt&\leq&\int_0^T\int_0^t\int_{\mathbb R^3}e^{-\nu(t-\tau
)|\xi|^2}|\xi|^2e^{a|\xi|^{1/\sigma}}|\widehat{u\otimes u}(\tau,\xi)|d\tau dtd\xi\\&\leq&\int_{\mathbb R^3}|\xi|^2e^{a|\xi|^{1/\sigma}}\left(\int_0^T
\int_0^te^{-\nu(t-\tau)|\xi|^2}|\widehat{u\otimes u}(\tau,\xi)|d\tau dt\right)d\xi.
\end{eqnarray*}
Integrating the function $e^{-\nu(t-\tau)|\xi|^2}$ twice with respect to $\tau \in[0,t]$ and $t\in[0,T]$, we get
\begin{eqnarray*}
\int_0^T\int_0^te^{-\nu(t-\tau)|\xi|^2}|\widehat{u\otimes u}(\tau,\xi)|d\tau dt&=&\int_0^T|\widehat{u\otimes u}(\tau,\xi)|\left(\left[\frac{-e^{-\nu(t
-\tau)|\xi|^2}}{\nu|\xi|^2}\right]_{\tau}^T\right)d\tau\\&\leq&\int_0^T|\widehat{u\otimes u}(\tau,\xi)|(\frac{1-e^{-\nu(T-\tau)|\xi|^2}}{\nu
|\xi|^2})d\tau.
\end{eqnarray*}
Then
\begin{eqnarray*}
\int_0^T\|\int_0^t e^{\nu(t-\tau)\Delta}div(u\otimes u)d\tau\|_{Z^{1}_{a,\sigma}}dt&\leq&\int_{\mathbb R^3}|\xi|^2e^{a|\xi|^{1/\sigma}}
\left(\int_0^T(\frac{1-e^{-\nu(T-\tau)|\xi|^2}}{\nu|\xi|^2})|\widehat{u\otimes u}(\tau,\xi)|d\tau\right)d\xi\\&\leq&\int_0^T\|u\otimes u\|_{Z^{0}
_{a,\sigma}}.
\end{eqnarray*}
Using the lemma \ref{lem1}, we will get the result.
\hfill $\square$\\
The proof of the first main result requires the following lemma.
\begin{lem}\label{lem4}
$$
\|u\otimes u\|_{Z^{0}_{a,\sigma}}\leq\|u\|_{Z^{-1}_{\frac{a}{\sqrt{\sigma}},\sigma}}\|u\|_{Z^{-1}_{a,\sigma}}^{\frac{1}{2}}\|\Delta u\|_{Z^{-1}_{a,\sigma}}^{\frac{1}{2}}.
$$
\end{lem}
\noindent{\it Proof lemma \ref{lem4}.}\\
It is easy to see that
$$
x^2e^{(\frac{a}{\sigma}-\frac{a}{\sqrt{\sigma}})x^{\frac{1}{\sigma}}}\leq c_{a,\sigma},\,\,\,\forall x\geq0.
$$
Then, for $x=|\xi|$
$$
|\xi|^2e^{(\frac{a}{\sigma}-\frac{a}{\sqrt{\sigma}})|\xi|^{\frac{1}{\sigma}}}\leq c_{a,\sigma}.
$$
This implies
$$
|\xi|e^{\frac{a}{\sigma}|\xi|^{\frac{1}{\sigma}}}\leq c_{a,\sigma}\frac{1}{|\xi|}e^{\frac{a}{\sqrt{\sigma}}|\xi|^{\frac{1}{\sigma}}}.
$$
Then
\begin{eqnarray*}
\|\Delta u\|_{Z^{-1}_{\frac{a}{\sigma},\sigma}}&=&\int_{\mathbb R^3}|\xi|e^{\frac{a}{\sigma}|\xi|^{\frac{1}{\sigma}}}|\hat{u}(\xi)|d\xi
\\&\leq&c\int_{\mathbb R^3}\frac{1}{|\xi|}e^{\frac{a}{\sqrt{\sigma}}|\xi|^{\frac{1}{\sigma}}}|\hat{u}(\xi)|d\xi\\&\leq&c\|u\|_{Z^{-1}_{\frac{a}{\sqrt{
\sigma}},\sigma}}.
\end{eqnarray*}
Using the previous computations and Cauchy-Schwartz inequality, we get
\begin{eqnarray*}
\|u\otimes u\|_{Z^{0}_{a,\sigma}}&=&\int_{\xi}e^{a|\xi|^{1/\sigma}}(\int_{\eta}|\hat{u}(\xi-\eta)||\hat{u}(\eta)|d\eta)d\xi\\&\leq&
c\|u\|_{Z^{0}_{\frac{a}{\sigma},\sigma}}\|u\|_{Z^{0}_{a,\sigma}}\\&\leq&c\|\Delta u\|_{Z^{-1}_{\frac{a}{\sigma},\sigma}} \|u\|_{Z^{0}_{a,\sigma}}\\&\leq&c\|u\|_{Z^{-1}_{\frac{a}{\sqrt{\sigma}},\sigma}}\|u\|_{Z^{0}_{a,\sigma}}
\\&\leq&c\|u\|_{Z^{-1}_{\frac{a}{\sqrt{\sigma}},\sigma}}\|u\|_{Z^{-1}_{a,\sigma}}^{\frac{1}{2}}\|\Delta u\|_{Z^{-1}_{a,\sigma}
}^{\frac{1}{2}}.
\end{eqnarray*}

\section{Well-posedness of $(NSE)$ in $Z^{-1}_{a,\sigma}(\mathbb{R}^3)$}
In the following theorem, we study the existence and uniqueness of the solution.
\begin{theorem}\label{theo5}
Let $u^0\in Z^{-1}_{a,\sigma}$. Then, there are a time $T>0$ and a unique solution $u\in{\mathcal C}([0,T],Z^{-1}_{a,\sigma}(\mathbb R^3))$ of $(NSE)$ such that $u\in{L^1}([0,T],Z^1_{a,\sigma}(\mathbb R^3))$.
\end{theorem}
\noindent{\it Proof theorem \ref{theo5}.}\\
{\bf(i)}Firstly, we wish to prove the existence.\\
The idea of the proof is to write the initial condition as a sum of higher and lower frequencies. For small frequencies, we will give a regular solution of the associated linear system to $(NSE)$. For the higher frequencies, we consider a partial differential equation very small to $(NSE)$ with small initial data in  $Z^{-1}_{a,\sigma}(\mathbb{R}^3)$ for which we can solve it by the Fixed Point Theorem.\\
$\bullet$ Let $r\in (0,\frac{1}{10})$.\\
$\bullet$ Let $N\in\mathbb{N}$, such that
$$
\int_{|\xi|>N}\frac{e^{a|\xi|^{1/\sigma}}}{|\xi|}|\hat{u^0}(\xi)|d\xi<\frac{r}{5}.
$$
Let's
$$
v^0=\mathcal{F}^{-1}(\mathbf{1_{\{|\xi|<N\}}}\hat{u^0}(\xi))
$$
and
$$
w^0=\mathcal{F}^{-1}(\mathbf{1_{\{|\xi|>N\}}}\hat{u^0}(\xi)).
$$
Clearly
\begin{eqnarray}\label{enq1}
\|w^0\|_{Z^{-1}_{a,\sigma}}<\frac{r}{5}.
\end{eqnarray}
Let $v=e^{\nu t\Delta}v^0$ the unique solution to
$$
\left\{
  \begin{array}{lll}
     \partial_t v
 -\nu\Delta v&=&\;\;0\\
    v(0,x) &=&v^0(x),
  \end{array}
\right.
\leqno{}
$$
We have
$$
\|v\|_{Z^{-1}_{a,\sigma}}\leq\|u^0\|_{Z^{-1}_{a,\sigma}},\,\,\forall t\geq0,
$$
and
\begin{eqnarray*}
\|v\|_{L^1_T(Z^{1}_{a,\sigma})}&=&\int_0^T\int_{\mathbb{R}^3}|\xi|e^{a|\xi|^{1/\sigma}}|\hat{v}(\xi)|d\xi dt\\&\leq&\int_0^T\int_
{\mathbb{R}^3}e^{-\nu t|\xi|^2}|\xi|e^{a|\xi|^{1/\sigma}}|\hat{u^0}(\xi)|d\xi dt\\&\leq&\int_{\mathbb{R}^3}(\int_0^Te^{-\nu t|\xi|^2}dt)
|\xi|e^{a|\xi|^{1/\sigma}}|\hat{u^0}(\xi)|d\xi\\&\leq&\frac{1}{\nu}\int_{\mathbb{R}^3}(1-e^{-\nu T|\xi|^2})|\xi|^{-1}e^{a|\xi|^{1/\sigma}} |\hat{u^0}(\xi)|d\xi.
\end{eqnarray*}
Using the Dominated Convergence Theorem, we get
\begin{eqnarray}\label{enq2}
\lim_{t\rightarrow0^+}\|v\|_{L^1_T(Z^{1}_{a,\sigma})}=0.
\end{eqnarray}
Let $\varepsilon>0$ such that
\begin{eqnarray*}
2\varepsilon\|u^0\|_{Z^{-1}_{a,\sigma}}<\frac{r}{5},
\end{eqnarray*}
\begin{eqnarray*}
\|u^0\|_{Z^{-1}_{a,\sigma}}+\varepsilon<\frac{1}{5},
\end{eqnarray*}
and
\begin{eqnarray*}
4(\varepsilon+2r\|u^0\|_{Z^{-1}_{a,\sigma}})\leq\frac{1}{2}.
\end{eqnarray*}
By \eqref{enq2}, there is a time $T=T(\varepsilon)>0$ such that
\begin{eqnarray*}
\|v\|_{L^1_T(Z^{1}_{a,\sigma})}<\varepsilon.
\end{eqnarray*}
Put $w=u-v$, clearly $w$ is the solution of the following system
$$
\left\{
  \begin{array}{lll}
     \partial_t w
 -\nu\Delta w+ (v+w).\nabla(v+w)  & =&\;\;-\nabla p\\
    w(0,x) &=&w^0(x) \;\;,
  \end{array}
\right.
\leqno{}
$$
The integral form of $w$ is as follows
$$
w=e^{\nu t\Delta}w^0-\int_0^t e^{\nu (t-\tau)\Delta}(v+w).\nabla(v+w)d\tau.
$$
To prove the existence of $w$, we put the following operator
$$
\psi(w)=e^{\nu t\Delta}w^0-\int_0^t e^{\nu (t-\tau)\Delta}(v+w).\nabla(v+w)d\tau.
$$
Now, we introduce the spaces $Z_T$ as follows
$$
Z_T=\mathcal C([0,T],Z^{-1}_{a,\sigma}(\mathbb R^3))\cap L^1([0,T],Z^{1}_{a,\sigma}(\mathbb R^3))
$$
with the norm
$$
\|f\|_{Z_T}=\|f\|_{L^{\infty}_T(Z^{-1}_{a,\sigma})}+\|f\|_{L^1_T(Z^{1}_{a,\sigma})}.
$$
Using lemmas \ref{lem2} and \ref{lem3}, we can prove $\psi(Z_T)\subset Z_T$.\\
$\bullet$ Also, denoted by $\mathbf{B}_r$ the subset of $Z_T$ defined by:
$$
\mathbf{B}_r=\{u\in Z_T;\|u\|_{L^{\infty}_T(Z^{-1}_{a,\sigma})}\leq r;\|u\|_{L^{1}_T(Z^{1}_{a,\sigma})}\leq r\}.
$$
$\bullet$ For $w\in\mathbf{B}_r$, we prove that $\psi(w)\subset\mathbf{B}_r$. In fact, we have
$$
\|\psi(w)(t)\|_{Z^{-1}_{a,\sigma}}\leq \sum_{k=0}^4 I_k,
$$
where
$$
I_0=\|e^{\nu t\Delta}w^0\|_{Z^{-1}_{a,\sigma}}
$$
$$
I_1=\int_0^t\|e^{\nu (t-\tau)\Delta}v\nabla v\|_{Z^{-1}_{a,\sigma}}d\tau
$$
$$
I_2=\int_0^t\|e^{\nu (t-\tau)\Delta}v\nabla w\|_{Z^{-1}_{a,\sigma}}d\tau
$$
$$
I_3=\int_0^t\|e^{\nu (t-\tau)\Delta}w\nabla v\|_{Z^{-1}_{a,\sigma}}d\tau
$$
$$
I_4=\int_0^t\|e^{\nu (t-\tau)\Delta}w\nabla w\|_{Z^{-1}_{a,\sigma}}d\tau.
$$
Using \eqref{enq1} the lemma \ref{lem2} and the fact that $w\in\mathbf{B}_r$, hence we get
$$
I_0\leq\frac{r}{5}
$$
\begin{eqnarray*}
I_1&\leq&\|v\|_{L^{\infty}_T(Z^{-1}_{a,\sigma})}\|v\|_{L^{1}_T(Z^{1}_{a,\sigma})}\\&\leq&2\varepsilon
\|u^0\|_{Z^{-1}_{a,\sigma}}\\&<&\frac{r}{5}
\end{eqnarray*}
\begin{eqnarray*}
I_2,I_3&\leq&\|v\|_{L^{1}_T(Z^{1}_{a,\sigma})}\|w\|_{L^{\infty}_T(Z^{-1}_{a,\sigma})}\\&+&\|v\|_{L^{\infty}_T(Z^{-1}_
{a,\sigma})}\|w\|_{L^{1}_T(Z^{1}_{a,\sigma})}\\&\leq&r(\|u^0\|_{Z^{-1}_{a,\sigma}}+\varepsilon)
\\&<&\frac{r}{5}
\end{eqnarray*}
\begin{eqnarray*}
I_4&\leq&\|w\|_{L^{\infty}_T(Z^{-1}_{a,\sigma})}\|w\|_{L^{1}_T(Z^{1}_{a,\sigma})}\\&\leq&2r^2\\&<&\frac{r}{5}.
\end{eqnarray*}
Then
\begin{eqnarray}\label{enq3}
\|\psi(w)(t)\|_{Z^{-1}_{a,\sigma}}\leq r.
\end{eqnarray}
Similarly, $$\|\psi(w)(t)\|_{L^1(Z^{1}_{a,\sigma})}\leq  \sum_{k=0}^4 J_k, $$ where
$$
J_0=\int_0^T\|e^{\nu t\Delta}w^0\|_{Z^{1}_{a,\sigma}}dt
$$
$$
J_1=\int_0^T\|\int_0^te^{\nu (t-\tau)\Delta}v\nabla v d\tau\|_{Z^{1}_{a,\sigma}}dt
$$
$$
J_2=\int_0^T\|\int_0^te^{\nu (t-\tau)\Delta}v\nabla w d\tau\|_{Z^{1}_{a,\sigma}}dt
$$
$$
J_3=\int_0^T\|\int_0^te^{\nu (t-\tau)\Delta}w\nabla v d\tau\|_{Z^{1}_{a,\sigma}}dt
$$
$$
J_4=\int_0^T\|\int_0^te^{\nu (t-\tau)\Delta}w\nabla w d\tau\|_{Z^{1}_{a,\sigma}}dt.
$$
Using lemmas \ref{lem3} and the fact that $w\in\mathbf{B}_r$, we get
$$
J_0\leq\frac{r}{5}
$$
\begin{eqnarray*}
J_1&\leq&2\|v\|_{L^{\infty}_T(Z^{-1}_{a,\sigma})}\|v\|_{L^{1}_T(Z^{1}_{a,\sigma})}\\&\leq&2\varepsilon
\|u^0\|_{Z^{-1}_{a,\sigma}}\\&<&\frac{r}{5}
\end{eqnarray*}
\begin{eqnarray*}
J_2,J_3&\leq&\|v\|_{L^{1}_T(Z^{1}_{a,\sigma})}\|w\|_{L^{\infty}_T(Z^{-1}_{a,\sigma})}\\&+&\|v\|_{L^{\infty}_T(Z^{-1}_
{a,\sigma})}\|w\|_{L^{1}_T(Z^{1}_{a,\sigma})}\\&\leq&r(\|u^0\|_{Z^{-1}_{a,\sigma}}+\varepsilon)
\\&<&\frac{r}{5}
\end{eqnarray*}
\begin{eqnarray*}
J_4&\leq&2\|w\|_{L^{\infty}_T(Z^{-1}_{a,\sigma})}\|w\|_{L^{1}_T(Z^{1}_{a,\sigma})}\\&\leq&2r^2\\&<&\frac{r}{5}.
\end{eqnarray*}
Then
\begin{eqnarray}\label{enq4}
\|\psi(w)(t)\|_{L^1(Z^{1}_{a,\sigma})}\leq r.
\end{eqnarray}
Combining \eqref{enq3} and \eqref{enq4}, we get  $\psi(w)\subset\mathbf{B}_r$ and we can deduce
\begin{eqnarray}\label{enq5}
\psi(\mathbf{B}_r)\subset\mathbf{B}_r.
\end{eqnarray}
$\bullet$ Proof of the following estimate
$$
\|\psi(w_2)-\psi(w_1)\|_{Z_T}\leq \frac{1}{2}\|w_2-w_1\||_{Z_T},\,\,\,w_1,w_2\in\mathbf{B}_r.
$$
In fact, we have
\begin{eqnarray*}
\psi(w_2)-\psi(w_1)&=&-\int_0^t e^{\nu(t-\tau)\Delta}((v+w_2)\nabla(v+w_2)-(v+w_1)\nabla(v+w_1))d\tau\\&=&-\int_0^t e^{\nu(t-\tau)\Delta}((v+w_2)\nabla
(w_2-w_1)+(w_2-w_1)\nabla(v+w_1))d\tau
\end{eqnarray*}
and
$$
\|\psi(w_2)-\psi(w_1)\|_{Z^{-1}_{a,\sigma}}\leq K_1+K_2,
$$
with
$$
K_1=\|\int_0^te^{\nu(t-\tau)\Delta}(v+w_2)\nabla(w_2-w_1)d\tau\|_{Z^{-1}_{a,\sigma}},
$$
$$
K_2=\|\int_0^te^{\nu(t-\tau)\Delta}(w_2-w_1)\nabla(v+w_1)d\tau\|_{Z^{-1}_{a,\sigma}}.
$$
Using lemma \ref{lem2}, we can deduce
\begin{eqnarray*}
K_1&\leq&\|v+w_2\|_{Z^{-1}_{a,\sigma}}\|w_2-w_1\||_{Z^{1}_{a,\sigma}}+\|v+w_2\|_{Z^{1}_{a,\sigma}}\|w_2-w_1\||
_{Z^{-1}_{a,\sigma}}\\&\leq&(\|v\|_{Z^{-1}_{a,\sigma}}+\|w_2\|_{Z^{-1}_{a,\sigma}})\|w_2-w_1\||_{Z^{1}_{a,\sigma
}}\\&+&(\|v\|_{Z^{1}_{a,\sigma}}+\|w_2\|_{Z^{1}_{a,\sigma}})\|w_2-w_1\||_{Z^{-1}_{a,\sigma}}\\&\leq&
(\varepsilon+2r+\|u^0\|_{Z^{-1}_{a,\sigma}})\|w_2-w_1\||_{Z_T}.
\end{eqnarray*}
Similarly, we get
$$
K_2\leq(\varepsilon+2r+\|u^0\|_{Z^{-1}_{a,\sigma}})\|w_2-w_1\||_{Z_T}.
$$
Then
\begin{eqnarray}\label{enq6}
\|\psi(w_2)-\psi(w_1)\|_{L^{\infty}_T(Z^{-1}_{a,\sigma})}\leq2(\varepsilon+2r+\|u^0\|_{Z^{-1}_{a,\sigma}})\|w_2-w_1\||_{Z_T}.
\end{eqnarray}
Therefore, we have
$$
\|\psi(w_2)-\psi(w_1)\|_{L^1(Z^{1}_{a,\sigma})}\leq K_3+K_4,
$$
with
$$
K_3=\int_0^T\|\int_0^te^{\nu(t-\tau)\Delta}(v+w_2)\nabla(w_2-w_1)d\tau\|_{Z^{1}_{a,\sigma}}dt,
$$
$$
K_4=\int_0^T\|\int_0^te^{\nu(t-\tau)\Delta}(w_2-w_1)\nabla(v+w_1)d\tau\|_{Z^{1}_{a,\sigma}}dt.
$$
Using lemma \ref{lem3}, then we can deduce
\begin{eqnarray*}
K_3&\leq&\|v+w_2\|_{L^1_T(Z^{1}_{a,\sigma})}\|w_2-w_1\||_{L^{\infty}_T(Z^{-1}_{a,\sigma})}\\&+&\|v+w_2\|_{L^{\infty}(Z^{-1}_{a,\sigma}
)}\|w_2-w_1\||_{L^1_T(Z^{-1}_{a,\sigma})}\\&\leq&(\varepsilon+2r+\|u^0\|_{Z^{-1}_{a,\sigma}})\|w_2-w_1\||_{Z_T}.
\end{eqnarray*}
Similarly, we get
$$
K_4\leq(\varepsilon+2r+\|u^0\|_{Z^{-1}_{a,\sigma}})\|w_2-w_1\||_{Z_T}.
$$
Then
\begin{eqnarray}\label{enq7}
\|\psi(w_2)-\psi(w_1)\|_{L^1_T(Z^{1}_{a,\sigma})}\leq2(\varepsilon+2r+\|u^0\|_{Z^{-1}_{a,\sigma}})\|w_2-w_1\||_{Z_T}.
\end{eqnarray}
By \eqref{enq6} and \eqref{enq7}, we obtain
$$
\|\psi(w_2)-\psi(w_1)\|_{Z_T}\leq4(\varepsilon+2r+\|u^0\|_{Z^{-1}_{a,\sigma}})\|w_2-w_1\||_{Z_T}.
$$
This implies
\begin{eqnarray}\label{enq8}
\|\psi(w_2)-\psi(w_1)\|_{Z_T}\leq\frac{1}{2}\|w_2-w_1\||_{Z_T}.
\end{eqnarray}
So, combining \eqref{enq5} and \eqref{enq8} and the Fixed Point Theorem, there is a unique $w\in\mathbf{B}_r$ such that $u=v+w$ is the solution of $(NSE)$ with $
u\in Z_T(\mathbb{R}^3)$.\\

{\bf(ii)} Secondly, we want to prove the uniqueness.\\
Let $u_1,u_2\in{\mathcal C}([0,T],Z^{-1}_{a,\sigma}(\mathbb R^3))\cap {L^1}([0,T],Z^1_{a,\sigma}(\mathbb R^3))$ of $(NSE)$ such that $ u_1(0)=u_2(0)$.
Put $\delta=U_1-U_2$. We have
\begin{eqnarray}\label{enq9}
 \partial_t\delta-\nu\Delta \delta+ u_1.\nabla \delta+\delta.\nabla u_2=-\nabla(p_1-p_2).
\end{eqnarray}
Then
$$
\partial_t\hat{\delta}+\nu|\xi|^2 \hat{\delta}+\widehat{ (u_1.\nabla \delta)}+\widehat{(\delta.\nabla u_2)}=0.
$$
Multiplying the previous equation by $\overline{\hat{\delta}}$, we get
\begin{eqnarray}\label{enq10}
\partial_t\hat{\delta}.\overline{\hat{\delta}}+\nu|\xi|^2 \hat{\delta}.\overline{\hat{\delta}}+\widehat{ (u_1.\nabla \delta)}.\overline{\hat{\delta}}
+\widehat{(\delta.\nabla u_2)}.\overline{\hat{\delta}}=0.
\end{eqnarray}
From Eq (\ref{enq9}) we have
$$
\partial_t\overline{\hat{\delta}}+\nu|\xi|^2 \overline{\hat{\delta}}+\overline{\widehat{ (u_1.\nabla \delta)}}+\overline{\widehat{(\delta.\nabla u_2)}}=0.
$$
Multiplying this equation by $\hat{\delta}$, we get
\begin{eqnarray}\label{enq11}
\partial_t\overline{\hat{\delta}}.\hat{\delta}+\nu|\xi|^2 \overline{\hat{\delta}}.\hat{\delta}+\overline{\widehat{ (u_1.\nabla \delta)}}.\hat{\delta}+\overline{\widehat{(\delta.\nabla u_2)}}.\hat{\delta}=0.
\end{eqnarray}
By summing (\ref{enq10}) and (\ref{enq11}), we get
$$
\partial_t|\hat{\delta}|^2+2\nu|\xi|^2|\hat{\delta}|^2+2Re(\widehat{ (u_1.\nabla \delta)}.\overline{\hat{\delta}})+2Re(\widehat{(\delta.\nabla u_2)}.\overline{\hat{\delta}})=0,
$$
and
$$
\partial_t|\hat{\delta}|^2+2\nu|\xi|^2|\hat{\delta}|^2\leq2|\widehat{ (u_1.\nabla \delta)}||\overline{\hat{\delta}}|+2|\widehat{(\delta.\nabla u_2)}||\overline{\hat{\delta}}|.
$$
Let $\varepsilon>0$, thereby we have
$$
\partial_t|\hat{\delta}|^2=\partial_t(|\hat{\delta}|^2+{\varepsilon}^2)=2\sqrt{|\hat{\delta}|^2+{\varepsilon}^2}.\partial_t\sqrt{|\hat{\delta}|^2+
{\varepsilon}^2}
$$
then
\begin{eqnarray*}
2\partial_t\sqrt{|\hat{\delta}|^2+{\varepsilon}^2}+2\nu|\xi|^2\frac{|\hat{\delta}|^2}{\sqrt{|\hat{\delta}|^2+{\varepsilon}^2}}&\leq&2|\widehat{ (u_1.\nabla \delta)}|\frac{|\hat{\delta}|}{\sqrt{|\hat{\delta}|^2+{\varepsilon}^2}}+2|\widehat{(\delta.\nabla u_2)}|\frac{|\hat{\delta}|}{\sqrt{|\hat{\delta}|^2+ {\varepsilon}^2}}\\&\leq&2|\widehat{ (u_1.\nabla\delta)}|+2|\widehat{(\delta.\nabla u_2)}|.
\end{eqnarray*}
By integrating with respect to time
$$
\sqrt{|\hat{\delta}|^2+{\varepsilon}^2}+\nu\int_0^t|\xi|^2\frac{|\hat{\delta}|^2}{\sqrt{|\hat{\delta}|^2+{\varepsilon}^2}}\leq\int_0^t|\widehat{ (u_1.\nabla\delta)}|d\tau+\int_0^t|\widehat{(\delta.\nabla u_2)}|d\tau.
$$
Letting $\varepsilon\rightarrow0$, we get
$$
|\hat{\delta}|+\nu\int_0^t|\xi|^2|\hat{\delta}|d\tau\leq\int_0^t|\widehat{ (u_1.\nabla\delta)}|d\tau+\int_0^t|\widehat{(\delta.\nabla u_2)}|d\tau.
$$
Multiplying by $\frac{e^{a|\xi|^{\frac{1}{\sigma}}}}{|\xi|}$ and integrating with respect to $\xi$ , thereafter we get
\begin{eqnarray*}
\|\delta\|_{Z^{-1}_{a,\sigma}}+\nu\int_0^t\|\Delta\delta\|_{Z^{-1}_{a,\sigma}}d\tau&\leq&\int_0^t\| u_1.\nabla\delta\|_{Z^{-1}_{a,\sigma}}d\tau+\int_0^t\|\delta.\nabla u_2\|_{Z^{-1}_{a,\sigma}}d\tau\\&\leq&\int_0^t\|
\delta u_1\|_{Z^{0}_{a,\sigma}}d\tau+\int_0^t\|u_2\delta\|_{Z^{0}_{a,\sigma}}d\tau.
\end{eqnarray*}
Using the elementary inequality $xy\leq\frac{x^2}{2}+\frac{y^2}{2}$, we get
\begin{eqnarray*}
\|\delta u_1\|_{Z^{0}_{a,\sigma}}&\leq&\|\delta\|_{Z^{0}_{a,\sigma}}\|u_1\|_{Z^{0}_{a,\sigma}}\\&\leq&\|\delta\|_{Z^{-1}_{a,\sigma}}^{\frac{1}{2}}\|\Delta\delta\|_
{Z^{-1}_{a,\sigma}
}^{\frac{1}{2}}\|u_1\|_{Z^{-1}_{a,\sigma}}^{\frac{1}{2}}\|\Delta u_1\|_{Z^{-1}_{a,\sigma}}^{\frac{1}{2}}\\&\leq&\frac{2}{\nu}
\|\delta\|_{Z^{-1}_{a,\sigma}}\|u_1\|_{Z^{-1}_{a,\sigma}}\|\Delta u_1\|_{Z^{-1}_{a,\sigma}}+\frac{\nu}{2}
\|\Delta\delta\|_{Z^{-1}_{a,\sigma}}.
\end{eqnarray*}
Similarly,
$$
\|u_2\delta\|_{Z^{0}_{a,\sigma}}\leq\frac{2}{\nu}\|\delta\|_{Z^{-1}_{a,\sigma}}\|u_2\|_{Z^{-1}_{a,\sigma}}\|\Delta u_2\|_{Z^{-1}_{a,\sigma}}+\frac{\nu}{2}\|\Delta\delta\|_{Z^{-1}_{a,\sigma}}.
$$
Then
\begin{eqnarray*}
\|\delta\|_{Z^{-1}_{a,\sigma}}&\leq&\frac{2}{\nu}\int_0^t\|\delta\|_{Z^{-1}_{a,\sigma}}\|u_1\|_{Z^{-1}_{a,\sigma}}\|\Delta u_1\|_{Z^{-1}_{a,\sigma}}d\tau\\&+&\frac{2}{\nu}\int_0^t\|\delta\|_{Z^{-1}_{a,\sigma}}\|u_2\|_{Z^{-1}_{a,\sigma}}\|\Delta u_2\|_{Z^{-1}_{a,\sigma}}d\tau.
\end{eqnarray*}
Using Gronwall lemma and the fact $(t\mapsto\|u_1\|_{Z^{-1}_{a,\sigma}}\|\Delta u_1\|_{Z^{-1}_{a,\sigma}}) \in L^1([0,T])$,
$(t\mapsto\|u_2\|_{Z^{-1}_{a,\sigma}}\|\Delta u_2\|_{Z^{-1}_{a,\sigma}}) \in L^1([0,T])$, we can deduce that $\delta=0$ in $[0,T]$
which gives the uniqueness.
\hfill $\square$\\

In the following, we prove a global existence if the initial condition is small in the Lei-Lin-Gevrey spaces.
\begin{theorem}\label{theo6}
Let $u^0\in Z^{-1}_{a,\sigma}(\mathbb R^3)$ such that $\|u^0\|_{Z^{-1}_{a,\sigma}}<\nu$. Then, there exists a unique global solution $u\in{\mathcal C}(\mathbb R^+,Z^{-1}_{a,\sigma}(\mathbb R^3))\cap {L^1}(\mathbb R^+,Z^1_{a,\sigma}(\mathbb R^3))$ of $(NSE)$ such that
$$
\|u(t)\|_{Z^{-1}_{a,\sigma}}+(\frac{\nu-\|u^0\|_{Z^{-1}_{a,\sigma}}}{2})\int_0^t\|\Delta u\|_{Z^{-1}_{a,\sigma}}d\tau\leq\|u^0\|_{Z^{-1}_{a,\sigma}}.
$$
\end{theorem}
\noindent{\it Proof theorem \ref{theo6}.}\\
From theorem \ref{theo3}, if $u^0\in Z^{-1}_{a,\sigma}(\mathbb R^3)$, we have a local existence
$$
u\in{L^{\infty}_T(Z^{-1}_{a,\sigma}(\mathbb R^3))}\cap {L^1_T(Z^1_{a,\sigma}(\mathbb R^3))}.
$$
Assume that $\|u^0\|_{Z^{-1}_{a,\sigma}}<\nu$ and $u\in{\mathcal C}([0,T^*),Z^{-1}_{a,\sigma}(\mathbb R^3))\cap {L^1_{loc}}([0,T^*),Z^1_{a,\sigma}(\mathbb R^3))$ is the maximal solution of $(NSE)$. We have
$$
\partial_t\|u(t)\|_{Z^{-1}_{a,\sigma}}+\nu\|\Delta u\|_{Z^{-1}_{a,\sigma}}\leq\|{\rm div}(u\otimes u)\|_{Z^{-1}_{a,\sigma}}.
$$
Integrating over $(0,t)$ we get
\begin{eqnarray}
\nonumber\|u(t)\|_{Z^{-1}_{a,\sigma}}+\nu\int_0^t\|\Delta u\|_{Z^{-1}_{a,\sigma}}d\tau&\leq&\|u^0\|_{Z^{-1}_{a,\sigma}}+
\int_0^t\|u\otimes u\|_{Z^{0}_{a,\sigma}}d\tau\\
&\leq&\|u^0\|_{Z^{-1}_{a,\sigma}}+\int_0^t\|u\|_{Z^{-1}_{a,\sigma}
}\|\Delta u\|_{Z^{-1}_{a,\sigma}}d\tau.\label{enq12}
\end{eqnarray}
Therefore, for $T_*=\sup \{t\in[0,T^*)\,\,/\,\|u(t)\|_{Z^{-1}_{a,\sigma}}<\alpha\}$, where $\alpha=\frac{\nu+\|u^0\|_{Z^{-1}_{a,\sigma}}}{2}$.\\
Take $t\in[0,T_*)$. Then we have
$$
\|u(t)\|_{Z^{-1}_{a,\sigma}}+\nu\int_0^t\|\Delta u\|_{Z^{-1}_{a,\sigma}}d\tau\leq\|u^0\|_{Z^{-1}_{a,\sigma}}
+\alpha\int_0^t\|\Delta u\|_{Z^{-1}_{a,\sigma}}d\tau.
$$
This implies
\begin{eqnarray*}
\|u(t)\|_{Z^{-1}_{a,\sigma}}+(\nu-\alpha)\int_0^t\|\Delta u\|_{Z^{-1}_{a,\sigma}}d\tau&\leq&\|u^0\|_{Z^{-1}_{a,\sigma}}\\&<&\alpha.
\end{eqnarray*}
Then $T_*=T^*$. Particularly if $T<T^*$, we have
$$
\|u(T)\|_{Z^{-1}_{a,\sigma}}+(\nu-\alpha)\int_0^T\|\Delta u\|_{Z^{-1}_{a,\sigma}}d\tau\leq\|u^0\|_{Z^{-1}_{a,\sigma}}.
$$
Therefore, $T^*=\infty$
\hfill $\square$

\section{Global solution}
In this section, we prove the first main theorem \ref{theo3}.\\
Let $u\in{\mathcal C}([0,T^*_{a,\sigma}),Z^{-1}_{a,\sigma}(\mathbb R^3))\cap {L^1_{loc}}([0,T^*_{a,\sigma}),Z^1_{a,\sigma}(\mathbb R^3))$ be the maximal solution of $(NSE)$, such that $\|u^0\|_{\mathcal{X}^{-1}}<\nu$.\\
Therefore, we have
\begin{eqnarray*}
\|u(t)\|_{Z^{-1}_{a,\sigma}}+\nu\int_0^t\|\Delta u\|_{Z^{-1}_{a,\sigma}}d\tau&\leq&\|u^0\|_{Z^{-1}_{a,\sigma}}+
\int_0^t\|div(u\otimes u)\|_{Z^{-1}_{a,\sigma}}d\tau\\&\leq&\|u^0\|_{Z^{-1}_{a,\sigma}}+\int_0^t\|u\otimes u\|_{Z^{0}_{a,\sigma}}d\tau.
\end{eqnarray*}
Using the lemma \ref{lem4} and the inequality $xy\leq\frac{x^2}{2}+\frac{y^2}{2}$, thus we get
\begin{eqnarray*}
\|u(t)\|_{Z^{-1}_{a,\sigma}}+\nu\int_0^t\|\Delta u\|_{Z^{-1}_{a,\sigma}}d\tau&\leq&\|u^0\|_{Z^{-1}_{a,\sigma}}+
c\int_0^t\|u\|_{Z^{-1}_{\frac{a}{\sqrt{\sigma}},\sigma}}\|u\|_{Z^{-1}_{a,\sigma}}^{\frac{1}{2}}\|\Delta u\|_{Z^{-1}_{a,\sigma}}^{\frac{1}{2}}d\tau\\&\leq&\|u^0\|_{Z^{-1}_{a,\sigma}}+c\int_0^t(\|u\|_{Z^{-1}_{\frac{a}{\sqrt{\sigma}},\sigma}}^2\|u\|_{Z^{-1} _{a,\sigma}}+\frac{\nu}{2}\|\Delta u\|_{Z^{-1}_{a,\sigma}})d\tau.
\end{eqnarray*}
This implies that
$$
\|u(t)\|_{Z^{-1}_{a,\sigma}}+\frac{\nu}{2}\int_0^t\|\Delta u\|_{Z^{-1}_{a,\sigma}}d\tau\leq\|u^0\|_{Z^{-1}_{a,\sigma}}+c\int_0^t\|u\|_{Z^{-1}_{\frac{a}{\sqrt{\sigma}},\sigma}}^2\|u\|_{Z^{-1}_{a,\sigma}}d\tau.
$$
By the Gronwall lemma, we get
$$
\|u(t)\|_{Z^{-1}_{a,\sigma}}\leq\|u^0\|_{Z^{-1}_{a,\sigma}}\exp(c\int_0^t\|u\|_{Z^{-1}_{\frac{a}{\sqrt{\sigma}},\sigma}}^2d\tau).
$$
Then
\begin{eqnarray*}
\|u(t)\|_{Z^{-1}_{a,\sigma}}+\nu\int_0^t\|\Delta u\|_{Z^{-1}_{a,\sigma}}d\tau&\leq&\|u^0\|_{Z^{-1}_{a,\sigma}}
\\&+&c\int_0^t\|u\|_{Z^{-1}_{\frac{a}{\sqrt{\sigma}},\sigma}}^2\|u^0\|_{Z^{-1}_{a,\sigma}}\exp(c\int_0^s\|u\|_{Z^{-1}_{\frac{a}{\sqrt{\sigma}},\sigma}
}^2)\\&\leq&\|u^0\|_{Z^{-1}_{a,\sigma}}(1+c\int_0^t\|u\|_{Z^{-1}_{\frac{a}{\sqrt{\sigma}},\sigma}}^2\exp(c\int_0^s
\|u\|_{Z^{-1}_{\frac{a}{\sqrt{\sigma}},\sigma}}^2))\\&\leq&\|u^0\|_{Z^{-1}_{a,\sigma}}\exp(c\int_0^t\|u\|_{Z^{-1}_{\frac{a}{\sqrt{\sigma}},\sigma}}^2)d\tau.
\end{eqnarray*}
Assumed that $T^*_{a,\sigma}<\infty$, by the previous inequality $\int_0^{T^*_{a,\sigma}}\|\Delta u\|_{Z^{-1}_{a,\sigma}}d\tau=\infty$. This implies that
$$
\int_0^{T^*_{a,\sigma}}\|u\|_{Z^{-1}_{\frac{a}{\sqrt{\sigma}},\sigma}}^2d\tau=\infty.
$$
As $Z^{-1}_{a,\sigma}(\mathbb R^3)\hookrightarrow Z^{-1}_{\frac{a}{\sqrt{\sigma}},\sigma}(\mathbb R^3)$. Then $T^*_{a,\sigma}=T^*_{\frac{a}{\sqrt{\sigma}},\sigma}$. Thus
\begin{eqnarray}\label{enq13}
T^*_{a,\sigma}=T^*_{\frac{a}{\sqrt{\sigma}},\sigma}=...=T^*_{\frac{a}{{\sigma}^{\frac{n}{2}}},\sigma},\,\,\,\forall n\in\mathbf{N}.
\end{eqnarray}
Therefore, from the dominated convergence theorem
$$
\lim_{n\rightarrow\infty}\|u^0\|_{Z^{-1}_{\frac{a}{{\sigma}^{\frac{n}{2}}},\sigma}}=\|u^0\|_{\mathcal{X}^{-1}}<\nu.
$$
Then, there exists $n_0\in\mathbf{N}$ such that
$$
\|u^0\|_{Z^{-1}_{\frac{a}{{\sigma}^{\frac{n}{2}}},\sigma}}<\nu,\,\,\,\,\forall n\geq n_0.
$$
Applying theorem \ref{theo4}, so we have $\forall n\geq n_0$
\begin{eqnarray}\label{enq14}
u\in{\mathcal C}(\mathbb R^+,Z^{-1}_{\frac{a}{{\sigma}^{\frac{n_0}{2}}},\sigma}).
\end{eqnarray}
Using the inequalities (\ref{enq13})-(\ref{enq14}) and for $n=n_0$, we obtain $ T^*_{a,\sigma}=T^*_{\frac{a}{{\sigma}^{\frac{n}{2}}},\sigma}=\infty$. This is absurd, so $ T^*_{a,\sigma}=\infty$.
\hfill $\square$

\section{Long time decay for the global solution}
In this section, we prove the second main theorem \ref{theo4}.\\
Let $u\in{\mathcal C}(\mathbb R^+,Z^{-1}_{a,\sigma}(\mathbb R^3))$. As $Z^{-1}_{a,\sigma}(\mathbb R^3)\hookrightarrow \mathcal{X}^{-1}(\mathbb R^3)$.
Then  $u\in{\mathcal C}(\mathbb R^+,\mathcal{X}^{-1}(\mathbb R^3))$.\\
For the results of Hantaek Bae (see \cite{HB}). There exist $t_0>0$ and $\alpha>0$ such that
\begin{eqnarray}\label{enq15}
\|e^{\alpha|D|}u(t)\|_{\mathcal{X}^{-1}(\mathbb R^3)}\leq c_0,\,\,\forall t\geq t_0,
\end{eqnarray}
where $t_0=\varphi(t)=\sqrt{t-t_0}$.\\
Therefore, let $a>0$ and $\beta>0$. Then, there exists $c_1>0$ such that
\begin{eqnarray}\label{enq16}
ax^{\frac{1}{\sigma}}\leq c_1+\beta x,\,\,x\geq0.
\end{eqnarray}
Take $\beta=\frac{\alpha}{2}$ and using the inequalities (\ref{enq15})-(\ref{enq16}) and the Cauchy-Schwartz inequality, so we obtain
\begin{eqnarray*}
\|u(t)\|_{Z^{-1}_{a,\sigma}}&=&\int_{\mathbb R^3}\frac{e^{a|\xi|^{1/\sigma}}}{|\xi|}|\hat{u}(\xi)|d\xi\\&\leq&\int_{\mathbb R^3}\frac {e^{c_1+\beta|\xi|}}{|\xi|}|\hat{u}(\xi)|d\xi\\&\leq&e^{c_1}\int_{\mathbb R^3}\frac{e^{\beta|\xi|}}{|\xi|}|\hat{u}(\xi)|d\xi\\&\leq&e^{c_1} \|e^{\alpha|D|}u(t)\|_{\mathcal{X}^{-1}}\|u\|^{\frac{1}{2}}_{\mathcal{X}^{-1}}\\&\leq&c_0 e^{c_1}\|u\|^{\frac{1}{2}}_{\mathcal{X}^{-1}}.
\end{eqnarray*}
Using theorem \ref{theo1}. So, $\lim_{t\rightarrow\infty}\|u(t)\|_{Z^{-1}_{a,\sigma}}=0$.
\hfill $\square$


\begin{thebibliography}{10}

\bibitem{MC1} M. Cannone, {\it Harmonic analysis tools for solving the incompressible Navier-Stokes equations\/}, Diterot Editeur, Paris, 1995.

\bibitem{MC2} M. Cannone, {\it Ondelettes,paraproduit et Navier-Stokes. Harmonic analysis tools for solving the incompressible Navier-Stokes equations\/},
 in: SFriedlander.D.Serra(Eds.), Handbook of Mathematical Fluid Dynamics, vol.3, Elsevier,2003.

\bibitem{JYC1}J-Y.Chemin, {\it Remarque sur l'existence global pour le systeme de Navier-Stokes incompressible\/}, SIAM J.Math.Anal.26 (2) (2009)
 599-624.

\bibitem{JYC2}J-Y.Chemin.I.Gallagher. {\it  Well-posedness and stability results for the Navier-Stokes in$\mathbb R^3$\/}, Ann, Inst.H.Poincare Anal Non Linear 26 (2)(2009) 599-624.

\bibitem{TK3}  H.Fujita,T.Kato, {\it On the Navier-Stokes initial value problem \/}, I.Arch.Ration.Mech.Anal.16 (1964) 269-315.

\bibitem{JB} J. Benameur, {\it Long Time Decay to the Lei-Lin solution of 3D Navier Stokes equation.\/} J.Math.Anal.Appl.(2015).

\bibitem{EH} E.Hopf. {\it Uber die Anfangswertaufgabe fur die hydrodinamischen Grundgleichungen\/}, Math.Nachr.4 (1951) 213-231.

\bibitem{TK1} T.Kato, {\it Quasi-Linear Equations of Evolution, With Application to Partial Differential Equations \/}, Lecture Notes in Math, vol.448,Sringer-Verlag,1975,pp. 25-70.

\bibitem{TK2} T.Kato. {\it $L^p$ -solution of the Navier Stokes in $\mathbb R^m$. With applications to weak solutions\/}, Math.Z.
187 (4)(1984) 471-480.

\bibitem{HKD} H.Koch.D.Tataru. {\it Well-posedness for the Navier Stokes equations \/}, Adv.Math.157(1)(2001) 22-35.

\bibitem{ZL} Z.Lei.F.Lin, {\it Global mild solutions of Navier Stokes equations  \/}, Comm.Pure Appl.Math.LXIV (2011) 1297 1304.

\bibitem{JL01} J.Leray. {\it Essai sur lr movement d'un liquide visqueux emplissant l'espace\/}, Acta Math.63 (1933) 22-25.

\bibitem{JL02} J.leray. {\it Sur le movement d'un liquide visqueux emplissant l'espace\/}, Acta Math.63 (1) (1934) 193-248.

\bibitem{FP} F.Planchon. {\it Global strong solutions in Sobolev or Lebesgue spaces to the incompressible Navier Stokes equations in
$\mathbb R^3 $\/},Ann.Inst.H.Poincare Anal.Non Lineaire 13 (3)(1996) 319 336.

\bibitem{HB} Hanteak. Bae. {\it Existence and Analyticity of Lie-Lin Solution to the Navier Stokes equations.\/}

\bibitem{ZZ} Z.Zhang,Z.Yin.  {\it Global well-posedness for the generalized Navier Stokes system.\/} arXiv:1306.3735v1 [Math.Ap]
17 june 2013.
\end{thebibliography}
\end{document}